\RequirePackage{ifpdf}
\ifpdf % We are running pdfTeX in pdf mode
\documentclass[pdftex]{sigma}
\else
\documentclass{sigma}
\fi

\begin{document}

\allowdisplaybreaks

\renewcommand{\PaperNumber}{030}

\FirstPageHeading

\ShortArticleName{Baker--Akhiezer Modules on Rational Varieties}

\ArticleName{Baker--Akhiezer Modules on Rational Varieties}

\Author{Irina A.~MELNIK~$^\dag$ and Andrey E.~MIRONOV~$^\ddag$}

\AuthorNameForHeading{I.A.~Melnik and A.E.~Mironov}

\Address{$^\dag$~Novosibirsk State University, 630090 Novosibirsk, Russia}
\EmailD{\href{mailto:sibirochka@ngs.ru}{sibirochka@ngs.ru}}

\Address{$^\ddag$~Sobolev Institute of Mathematics, 630090 Novosibirsk, Russia}
\EmailD{\href{mailto:mironov@math.nsc.ru}{mironov@math.nsc.ru}}
\URLaddressD{\url{http://math.nsc.ru/LBRT/d6/mironov/}}

\ArticleDates{Received January 05, 2010, in f\/inal form April 03, 2010;  Published online April 07, 2010}

\Abstract{The free Baker--Akhiezer modules on rational varieties obtained from
${\mathbb C}P^{1}\times{\mathbb C}P^{n-1}$ by identif\/ication of two
hypersurfaces  are constructed. The corollary of this construction
is the existence of embedding of meromorphic function ring with some
f\/ixed pole into the ring of matrix dif\/ferential operators in $n$
variables.}


\Classification{14H70; 35P30}

\section{Introduction}

The Baker--Akhiezer modules (BA-modules) over the ring of
dif\/ferential operators were introduced by Nakayashiki (see~\cite{N1,N2}). These modules are constructed on the basis of
spectral data which include an algebraic variety $X$ and some
additional objects. In the one-dimensional case the module elements
are the usual Baker--Akhiezer functions.

 The BA-module $M$ consists of functions
$\psi(x,P)$ which depend on $x \in {\mathbb C}^n$, where
$n=\dim_{\mathbb C} X$ and $P\in X$. If $x$ is f\/ixed, then the
function $\psi$ is the section of a bundle over $X$, and $\psi$ has
an essential singularity on  divisor $Y \subset X$. The elements
$\psi \in M$ have the following properties:
\begin{itemize}\itemsep=0pt
\item $\partial_{x_j} \psi \in M$ and $f(x) \psi \in M$, where $f(x)$
is an analytical function in a neighbourhood of a~f\/ixed point
$x_0$;
\item if
$\lambda$ is a meromorphic function with a pole on $Y$, then
$\lambda\psi\in M$.
\end{itemize}
These properties mean that $M$ is the module over the ring of
dif\/ferential operators $\mathcal{D}_n =
\mathcal{O}[\partial_{x_1},\dots,\partial_{x_n}]$, where
$\mathcal{O}$ is the ring of analytical functions in a neighbourhood
of $x_0,$ and over the ring $A_Y$ of meromorphic functions on $X$
with the pole on $Y$.

The f\/initely generated free BA-modules over $\mathcal{D}_n$ are of
the main interest  as in this case the construction allows to
build commutative rings of dif\/ferential operators. Let us choose
 the basis $\psi_1(x,P),\dots, \psi_{N}(x,P)$ in $M$. Let
$\Psi(x,P)$ denote vector-function
$(\psi_1(x,P),\dots,\psi_{N}(x,P))^{\top}$. Then for $\lambda\in
A_Y$ there is only one dif\/ferential operator $D(\lambda)$ with $N
\times N$-matrix coef\/f\/icients such that
\[
 D(\lambda)\Psi(x,P)=\lambda(P)\Psi(x,P).
\]
The operators $D(\lambda)$ and $D(\mu)$ obviously commutate with
each other for dif\/ferent $\lambda$ and $\mu \in A_Y$. Thus, the
considered construction makes it possible to obtain solutions of
nonlinear dif\/ferential equations which are equivalent to the
condition of commutation of dif\/ferential operators.

The following examples of the free BA-modules over $\mathcal{D}_n$
are known. In \cite{N1} and \cite{N2} it is shown that the
BA-modules on Abelian varieties are free under some restrictions on
spectral data. In~\cite{N2} it is also shown that the restriction of
the BA-module from three-dimensional Abelian variety to the shifted
theta-divisor remains free (over the ring of dif\/ferential operators
of two variables).

In \cite{M} it is shown\footnote{In the proof of freeness of the
BA-modules in  \cite{M} there is a gap, an additional proposition
is required. The complete proof will appear in the work of K.~Cho,
A.~Mironov, and A.~Nakayashiki ``Baker--Akhiezer module on the
intersection of shifted theta divisors'' (submitted to Publ. RIMS).}
that the restriction of BA-module from Abelian variety to the
complete intersection of the shifted theta-divisors remains free. On
these grounds we found the solutions of multidimensional analog of
 Kadomtsev--Petviashvili hierarchy
\begin{equation}
[\partial_{t_k}-A_k,\partial_{t_m}-A_m] = 0, \label{eq1}
\end{equation}
where $A_k$ and $A_m$ are matrix dif\/ferential operators of $n$
variables.

In \cite{R1} and \cite{R2} the suf\/f\/icient conditions for the
spectral data, which correspond to the free BA-modules  were found
(see Theorem 4.1 in \cite{R1} and Theorem 3.3 in \cite{R2}). It is
not clear how to f\/ind the algebraic varieties, satisfying these
conditions. In \cite{R1} and \cite{R2} there are two examples
satisfying these conditions. In
 \cite{R2}
suf\/f\/icient conditions for spectral data corresponding to the
solutions of the equation (\ref{eq1})  were ascertained; in the
paper a corresponding example was also demonstrated.

Note that in all the  examples above the construction is either
implicit (see \cite{R1} and \cite{R2}), or the solutions are
expressed in terms of theta-functions (see \cite{N1,N2},
and \cite{M}).

Those who wish to read more widely in the theory of commuting
operators of several variables and BA-modules can turn to~\cite{P}.

The aim of this paper is to demonstrate the construction of
Nakayashiki for the rational varie\-ties.  For the rational spectral
variety, the BA-module elements and coef\/f\/icients of commutating
dif\/ferential operators are expressed in terms of elementary
functions.

Our initial idea was to obtain BA-modules on rational varieties from
BA-modules of Nakaya\-shi\-ki by degenerating of Abelian varieties in
the same way as soliton solutions of KdV are obtained from
f\/inite-gap solutions by degenerating of smooth spectral curves to
spheres with double points. We considered many candidates of
rational varieties and, as a result, we found varieties $\Gamma$ and
$\Omega$ (see below) appropriate for our goals.

In the next section we describe the spectral data used in this paper
and formulate our main results. In Sections~\ref{section3} and~\ref{section4} we show that the
BA-modules on $\Gamma$ and $\Omega$ are free. In Section~\ref{section5} we
present explicit examples of commutating operators. In the Appendix
we show that on $\Gamma$ and~$\Omega$ there are structures of
algebraic varieties (it not follows directly from the def\/inition of
$\Gamma$ and~$\Omega$).

\section{Main results}\label{section2}

Let us f\/ix  $a_1, a_2, b_1, b_2 \in {\mathbb C}$ such that
$(a_{i},b_{i})\ne (0,0)$ and $(a_1:b_1)
\neq (a_2:b_2)$. Let us also f\/ix  nondegenerate linear map
$\mathcal{P}:{\mathbb C}^{n}\rightarrow{\mathbb C}^{n}$. This map
induce the map ${\mathbb C}P^{n-1}\rightarrow {\mathbb C}P^{n-1}$,
which we denote by the same symbol $\mathcal{P}$. Let $\Gamma$
denote the variety constructed from ${\mathbb C}P^{1}\times{\mathbb
C}P^{n-1}$ by identif\/ication of two
 hypersurfaces
\[
 p_{1}\times{\mathbb
 C}P^{n-1}\sim p_{2}\times{\mathbb C}P^{n-1}
\]
with the use of $\mathcal{P}$, where $p_{i}=(a_{i}:b_{i})$. Namely,
let us identify
\[
(a_{1}:b_{1},t)\sim(a_{2}:b_{2},\mathcal{P}(t)), \qquad
 t=(t_{1}:\dots:t_{n})\in{\mathbb C}P^{n-1}.
\]

Let $f(P)$ be  the following function on ${\mathbb C}^{n+2}$
\begin{equation} \label{eq2}
f(z_{1},z_{2},t_{1},\dots,t_{n})=\sum_{i=1}^{n} (\alpha_{i} z_{1}
t_{i}+\beta_{i} z_{2} t_{i}), \qquad \alpha_{i},\beta_{i}\in{\mathbb
C}.
\end{equation}
such that the following identity takes place
\begin{equation} \label{eq3}
f(a_{1},b_{1},t)-A f(a_{2},b_{2},\mathcal{P}(t))=0
\end{equation}
for f\/ixed $A\in{\mathbb C}^*$ and every $t=(t_{1},\dots,t_{n}).$
Identity (\ref{eq3}) gives the restriction on the choice of~$\alpha_{i}$,~$\beta_{i}$. According to (\ref{eq3}), the equation
\[
f(z_{1}:z_{2},t_{1}:\dots:t_{n})=0
\]
correctly def\/ines a hypersurface in $\Gamma$.

We denote eigenvalues and eigenvectors of $\mathcal{P}$ by
$\lambda_j$ and $\mathbf{w}_{j}$ respectively. Henceforth, we assume
that
\begin{equation} \label{eq4}
 \lambda_j \neq \lambda_k  \qquad \mbox{at} \quad  j \neq k,
\end{equation}
and $f(P)$ is chosen such that
\begin{equation} \label{eq5}
f(a_{1},b_{1},\mathbf{w}_{j})\neq 0 ,\qquad j = 1,\dots,n.
\end{equation}
We introduce $n$ functions on ${\mathbb C}^{n+2}$
\[
f_{i}(z_{1},z_{2},t_{1},\dots,t_{n})=\sum_{k=1}^{n} \left (
\alpha_{ik} z_{1} t_{k}+\beta_{ik} z_{2} t_{k} \right )
\]
such that $f_i$ satisf\/ies the identity:
\begin{equation} \label{eq6}
\frac{f_{i}(a_{1},b_{1},t)}{f(a_{1},b_{1},t)}-
\frac{f_{i}(a_{2},b_{2},\mathcal{P}(t))}{f(a_{2},b_{2},
\mathcal{P}(t))}-c_{i}=0, \qquad c_{i}\in{\mathbb C}
\end{equation}
for  every $t=(t_{1},\dots,t_{n}) \in {\mathbb C}^{n}$.  By (\ref{eq3}),
this identity is equivalent to
\begin{equation} \label{eq7}
f_{i}(a_{1},b_{1},t) - A f_{i}(a_{2},b_{2},\mathcal{P}(t))-c_{i}
f(a_{1},b_{1},t)=0.
\end{equation}
The dimension of the space of these functions is equal to $(n+1)$.
We choose $f_{1},\dots,f_{n}$ such that $f_{1},\dots,f_{n}$ and $f$
are linearly independent. Moreover, we choose parameters
$(\alpha,\beta)$, $(\alpha_i,\beta_i)$ in a general position, that
means that the parameters belong to some open everywhere dense
domain (more precisely, such that equation (\ref{eq16.6}) has no
multiple solutions).

Let us f\/ix $\Lambda \in {\mathbb C}$. Let
\[
M_{\Gamma}(k)=\left\{\psi(x,
P)=\frac{h(x_{1},\dots,x_{n},P)}{f^{k}(P)}\exp
\left(\sum_{j=1}^{n}\frac{f_{j}(P)}{f(P)}x_{j}\right)\right\},
\]
where
\begin{equation} \label{eq8}
\psi(x,a_{1}:b_{1},t) - \Lambda \psi(x,a_{2}:b_{2},\mathcal{P}(t)) =
0
\end{equation}
for $t=(t_{1}:\dots:t_{n})\in{\mathbb C}P^{n-1}$. Here
$P=(z_{1}:z_{2},t)\in {\mathbb C}P^{1}\times{\mathbb C}P^{n-1}$
and $h(x,P)$ has the form
\begin{equation} \label{eq9}
h(x,P)=\sum_{0\leq j \leq k,\, |\alpha|=k} h_{j
\alpha}(x) z_{1}^{j}z_{2}^{k-j}  t^{\alpha},
\end{equation}
where $\alpha=(\alpha_{1},\dots,\alpha_{n})$,
$t^{\alpha}=t_{1}^{\alpha_{1}} \cdot \dots \cdot
t_{n}^{\alpha_{n}}$.

According to (\ref{eq8}), if $\psi\in M_{\Gamma}(k)$, then
$\partial_{x_{j}}\psi\in M_{\Gamma}(k)$. Consequently, we have $n$
mapping
\[
\partial_{x_{j}} : \  M_{\Gamma}(k) \rightarrow M_{\Gamma}(k+1), \qquad
j=1,\dots,n.
\]
Thereby, on the set
\[
M_{\Gamma}=\bigcup_{k=1}^{\infty}M_{\Gamma}(k)
\]
the structure of the BA-module over the ring of dif\/ferential
operators $ \mathcal{D}_{n}=\mathcal{O}
\left[\partial_{x_{1}},\dots,\partial_{x_{n}}\right]$ is def\/ined.

\begin{theorem} \label{T1}
 $M_{\Gamma}$ is a free
$\mathcal{D}_{n}$-module of rank $n$ generated by $n$ functions
from $M_{\Gamma}(1)$.
\end{theorem}

\begin{corollary}  There is a ring embedding
\[
D: \ A_{f} \rightarrow {\rm Mat}(n,\mathcal{D}_{n})
\]
of the rings of meromorphic functions on $\Gamma$ with the poles on
hypersurface $f = 0$ into the ring of differential operators in
variables $x_{1},\dots,\,x_{n}$ with the matrix coefficients of
size $n\times n$.
\end{corollary}

For $n = 2$ it is possible to consider another way of identif\/ication
of two curves in ${\mathbb C}P^1\times {\mathbb C}P^1$. Let $\Omega$
denote the variety which comes out from ${\mathbb C}P^1\times
{\mathbb C}P^1$ by the identif\/ication of two lines
\[
 p_1\times {\mathbb C}P^1 \sim {\mathbb C}P^1\times p_2.
\]
Videlicet, we identify the following points:
\[
 (p_1,t)\sim (\mathcal{P}(t),p_2),
\]
where $p_i,   t \in{\mathbb C}P^1$. We assume that
$\mathcal{P}(p_2) \neq p_1$. Therefore, in an appropriate coordinate
system, the variety $\Omega$ has the form
\begin{equation} \label{eq10}
\Omega = {\mathbb C}P^1\times {\mathbb C}P^1 / \{(1:0,t_1:t_2)\sim
(t_1:t_2,0:1)\}.
\end{equation}
Indeed,  on ${\mathbb C}P^1\times {\mathbb C}P^1$ we make the
following change of coordinates:
\[
(z',w')=(z,\mathcal{P}(w)),
\]
where $(z,w)$, $(z',w')$ are the old and the new coordinates on
${\mathbb C}P^1\times {\mathbb C}P^1$ respectively. Then $\Omega$ is
obtained by the identif\/ication of the points:
\[
(p_1,\mathcal{P}(t))\sim (\mathcal{P}(t),\mathcal{P}(p_2)).
\]
Now on each of components ${\mathbb C}P^1\times {\mathbb C}P^1$ we
do the same change of coordinates such that the points $p_1$ and
$\mathcal{P}(p_2)$ in the new system have coordinates $(1:0)$ and
$(0:1)$ respectively. In the new coordinates the variety $\Omega$
has the form~(\ref{eq10}).

 Let $g$ denote the following function
\[
 g(z_1,z_2,w_1,w_2)=\alpha z_1w_1+\beta z_1w_2+\gamma z_2w_1+\delta
 z_2w_2, \qquad \alpha,\beta,\gamma,\delta\in{\mathbb C}
\]
such that for $(t_1,t_2)\in{\mathbb C}^2$  the following identity is
fulf\/illed
\begin{equation} \label{eq11}
 g(1,0,t_1,t_2)-B g(t_1,t_2,0,1)=0,\;
\end{equation}
where $ B\in {\mathbb C}^*$ is f\/ixed. We assume that
\[
g(0,1,0,1) \neq 0.
\]
We introduce two more functions
\[
 g_i(z_1,z_2,w_1,w_2)=\alpha_i z_1w_1+\beta_i z_1w_2+\gamma_i
 z_2w_1+\delta_i z_2w_2,\qquad \alpha_i,\beta_i,\gamma_i,\delta_i\in{\mathbb C},\quad
 i=1,2
\]
such that for $(t_1,t_2)\in{\mathbb C}^2$  the identity is fulf\/illed
\begin{equation} \label{eq12}
 \frac{g_i(1,0,t_1,t_2)}{g(1,0,t_1,t_2)}-
 \frac{g_i(t_1,t_2,0,1)}{g(t_1,t_2,0,1)}-c_i=0,\qquad
 c_i\in {\mathbb C}.
\end{equation}
The dimension of the space of such functions is equal to 3.
According to (\ref{eq11}), the identity (\ref{eq12}) is equivalent
to
\[
 g_i(1,0,t_1,t_2)-B g_i(t_1,t_2,0,1)-c_i  g(1,0,t_1,t_2)=0.
\]
Let us choose $g_1$ and $g_2$ in such a way that $g_1$, $g_2$ and $g$
are linearly independent and the under radical expression in
(\ref{eq21}) does not vanish (this always can be achieved by the
inf\/initesimal changes of $c_1$ and $c_2$). Let
\[
 G_1(P)=\frac{g_1(P)}{g(P)},\qquad G_2(P)=\frac{g_2(P)}{g(P)}.
\]
Let us f\/ix $\Lambda\in{\mathbb C}$.  By $M_{\Omega}(k)$ we denote
the set of functions of the form
\[
 M_{\Omega}(k)=\left\{\varphi=\frac{\tilde{h}(x,y,P)}{g^k(P)}
 \exp\left(x G_1(P)+y G_2(P) \right)\right\}
\]
for which the identity
\[
 \varphi(x,y,1:0,t_1:t_2)-\Lambda\varphi(x,y,t_1:t_2,0:1)=0,
\]
  is fulf\/illed, where $\tilde{h}$ is the function of the
form (\ref{eq9}).

Let
\[
 M_{\Omega}=\bigcup_{k=1}^{\infty}M_{\Omega}(k),
\]
$M_{\Omega}$ is the module over
$\mathcal{D}=\mathcal{O}[\partial_x,\partial_y]$.

\begin{theorem} \label{T2}
 $M_{\Omega}$ is a free $\mathcal{D}$-module of the rank~$2$ generated by two functions from $M_{\Omega}(1)$.
\end{theorem}

Let $A_g$ denote the ring of the meromorphic functions on $\Omega$
with the pole on the curve def\/ined by the equation $g(P)=0$.

\begin{corollary}  There is a ring embedding
\[
 D: \ A_g\rightarrow {\rm Mat}(2,\mathcal{D})
\]
of $A_g$ into the ring of $2\times 2$-matrix differential operators
in variables $x$ and $y$.
\end{corollary}

\begin{remark} To be more precise, the freeness of
$M_{\Gamma}$ and $M_{\Omega}$ is a corollary of the fact that
corresponding graded ${\rm gr}\mathcal{D}$-modules are free, where
the graduation is induced by the degree of the operators and the
order of poles respectively. Below we virtually prove the freeness
of graded modules.
\end{remark}

\section[Proof of Theorem 1]{Proof of Theorem \ref{T1}}\label{section3}

\subsection[Combinatorial structure of $M_{\Gamma}$]{Combinatorial structure of $\boldsymbol{M_{\Gamma}}$}

We f\/ind the rank of the module $M_{\Gamma}(k)$ over $\mathcal{O}$.
The dimension of the space of functions $\{h(x,P)\}$ (see (\ref{eq9})) is
equal $C_{k+1}^{1} C_{k+n-1}^{n-1}$ (for f\/ixed $x$). The condition
(\ref{eq8}) with the help of (\ref{eq3}) and (\ref{eq6}) can be written in the equivalent
form
\[
h(x,a_{1},b_{1},t)-h(x,a_{2},b_{2},\mathcal{P}(t)) \Lambda A^{k}
e^{-cx}=0,
\]
where $cx=\sum_{j=1}^{n}c_{j}x_{j}$. This equality means that the
coef\/f\/icients of the homogeneous polynomial in $t_{1},\dots,t_{n}$ of
the degree $k$, situated in the left part, are equal to~$0$. It
gives $C_{k+n-1}^{n-1}$ restrictions on the choice of coef\/f\/icients
of  $h(x, P)$. Thereby,
\begin{equation}\label{eq13}
 {\rm rank}_{\mathcal{O}}M_{\Gamma}(k) = (k+1)  C_{k+n-1}^{n-1}- C_{k+n-1}^{n-1}= k
 C_{k+n-1}^{n-1}.
\end{equation}
Let $\mathcal{D}_{n}^{k-1}$ be the dif\/ferential operators of order
 $k-1$ at most. We have
\begin{equation}\label{eq14}
{\rm rank}_{\mathcal{O}} \mathcal{D}_{n}^{k-1}=
C_{k+n-1}^{n}=\frac{(k+n-1)!}{(k-1)!n!}=
\frac{k(k+1) \cdots (k+n-1)}{n!}=\frac{k}{n} C_{k+n-1}^{n-1}.
\end{equation}
Comparing (\ref{eq13}) and (\ref{eq14}), we can expect that $M_{\Gamma}$ is a free
module of the rank~$n$ generated by~$n$ functions from
$M_{\Gamma}(1)$.

\subsection[Module $N$]{Module $\boldsymbol{N}$}

Let us choose $n$ functions $\psi_{1},\dots,\psi_{n}\in
M_{\Gamma}(1)$ independent over $\mathcal{O}$
\[
\psi_{k}(x,P)=\frac{h_{k}(x,P)}{f(P)}\exp\left(\sum_{j=1}^{n}
\frac{f_{j}(P)}{f(P)}x_{j}\right),
\]
where
\[
h_{k}(x,P)=\sum_{i=1}^{n} \left(h^{1}_{ki}(x)  z_{1}
t_{i}+h^{2}_{ki}(x)  z_{2} t_{i}\right).
\]
Consider the module $N$ over $\mathcal{D}_{n}$ generated by the
functions $\psi_{1},\dots,\psi_{n}$
\[
N=\left\{\sum_{i=1}^{n}d_{i} \psi_{i} \, |\,    d_{i}\in
\mathcal{D}_{n} \right\}.
\]
We show that the module $N$ is free (Lemma~\ref{L1}) and as a consequence
from the combinatorial calculation we get that the modules
$M_{\Gamma}$ and $N$ coincide (Lemma~\ref{L3}).

\begin{lemma} \label{L1}
 $N$ is a free $\mathcal{D}_{n}$-module of rank $n$.
\end{lemma}

\begin{proof} Suppose that the assertion is not true,
i.e.\ there are  dif\/ferential operators $d_{1}, \dots ,d_{n} \in
\mathcal{D}_{n}$ such that
\begin{equation}\label{eq15}
d_{1}\psi_{1}+\dots+d_{n}\psi_{n} = 0,
\end{equation}
where{\samepage
\[
d_{j}= \sum_{\alpha : |\alpha|\leq K}a_{j
\alpha}(x) \partial^{\alpha}_{x},
\]
$K$ is maximal order of operators $d_{j}$.}

Let us divide (\ref{eq15}) by $\exp\big(\sum\frac{f_{j}}{f}x_{j}\big)$,
multiply by $f^{K+1}$ and restrict the received equality on the
hyperspace $f(P)=0$. We receive the following equality (for the
compactness of the record the arguments are skipped)
\begin{gather}
h_{1} \big(a_{1 (K,0,\dots,0)}f_{1}^{K}+
a_{1 (K-1,1,0,\dots,0)}
f_{1}^{K-1}f_{2}+\dots+a_{1 (0,\dots,0,K)}
f_{n}^{K}\big)+\cdots\nonumber\\
\qquad{}+
h_{n} \big(a_{n(K,0,\dots,0)}f_{1}^{K}+ a_{n(K-1,1,0,\dots,0)}
f_{1}^{K-1}f_{2}+\dots+a_{n(0,\dots,0,K)}f_{n}^{K}\big)=0.\label{eq16}
\end{gather}
The hypersurfaces
\begin{equation} \label{eq16.5}
 f(P)=0, \qquad f_{j}(P)=0,\qquad j=1,\dots,n,\quad j\neq k.
\end{equation}
($f_k(P)=0$ is left out) have $n$ points of intersections. Indeed,
let us consider (\ref{eq16.5}) as linear equations in
$t_1,\dots,t_n$. These equations have nonzero solutions if the
determinant $\Delta_j,$ composed from coef\/f\/icients (coef\/f\/icients are
linear forms $\alpha_iz_1+\beta_iz_2$ or
$\alpha_{si}z_1+\beta_{si}z_2$), equals 0
\begin{equation}\label{eq16.6}
 \Delta_j=0.
\end{equation}
So, (\ref{eq16.6}) is a homogeneous equation in $z_1$, $z_2$ of degree~$n$, and by our assumption has no multiple solutions. By
$P_{i}^{k}$, $i=1,\dots,n$ we denote  the intersection points of
hypersurfaces~(\ref{eq16.5}).

Let us substitute  $P_{i}^{n}$ in (\ref{eq16}), divide all equalities by $f_{n}^{K}$ ($f_{n}(P_{i}^{n}) \neq
 0$, see Lemma~\ref{L4}). We obtain a linear system of equations
on coef\/f\/icients $a_{k\,(0,\dots,0,K)}$ of operators $d_k$
\begin{gather}
    h_{1}(x,P^{n}_{1})a_{1\,(0,\dots,0,K)}+\dots+h_{n}(x,P^{n}_{1})a_{n\,(0,\dots,0,K)}=0, \nonumber\\
    \cdots \cdots\cdots\cdots\cdots\cdots\cdots\cdots\cdots\cdots\cdots\cdots\cdots\cdots\cdots\cdots\cdots\nonumber\\
    h_{1}(x,P^{n}_{n})a_{1\,(0,\dots,0,K)}+\dots+h_{n}(x,P^{n}_{n})a_{n\,(0,\dots,0,K)}=0.  \label{eq17}
\end{gather}
We need
\begin{lemma}\label{L2}
The inequality holds
\[
\det\left(%
\begin{array}{ccc}
  h_{1}(x,P^{k}_{1}) & \dots & h_{n}(x,P^{k}_{1}) \\
  \vdots             & \ddots & \vdots \\
  h_{1}(x,P^{k}_{n}) & \dots & h_{n}(x,P^{k}_{n})
\end{array}
\right)\ne 0.
\]
\end{lemma}

The proof of this lemma is given in Subsection~\ref{section3.3}.

By Lemma \ref{L2}, the solution of the system (\ref{eq17}) is
\[
 a_{1(0,\dots,0,K)}=0, \quad \dots, \quad
 a_{n(0,\dots,0,K)}=0.
\]
Similarly, for $k=1,\dots,n-1$ from $\det\left(
h_{j}(x,P^{k}_{i})\right)\ne 0$ it follows that the coef\/f\/icients
$a_{j(0\dots K \dots 0)}$ must vanish ($K$ is on the $k$-th
place).

To show that the coef\/f\/icients  $a_{j\,(0,\dots0,m,l)}$, where
$m+l=K$, vanish we should restrict (\ref{eq16}) on $f(P)=0$,
$f_{1}(P)=0$, $\dots$, $f_{n-2}(P)=0$
\begin{gather*}
h_{1} \big(a_{1(0,\dots,0,K-1,1)}f_{n-1}^{K-1}f_{n}+
\dots+a_{1 (0,\dots,0,1,K-1)}f_{n-1}f_{n}^{K-1}\big)+\cdots
\\
\qquad{}+
h_{n} \big(a_{n(0,\dots,0,K-1,1)}f_{n-1}^{K-1}f_{n}+
\dots+a_{n(0,\dots,0,1,K-1)}f_{n-1}f_{n}^{K-1}\big)=0.
\end{gather*}
Since the product $f_{n-1}(P)f_{n}(P)$ is not identically zero on
this set, then it can be used for dividing. Substituting the
points $P=P^{n-1}_{i}$ and $P=P^{n}_{i}$ we obtain sets of
equations of the form (\ref{eq17}) on  coef\/f\/icients
$a_{j\,(0,\dots,1,K-1)}$ and $a_{j\,(0,\dots,K-1,1)}$. Using Lemma~\ref{L2}, we conclude that the matrices $\left( h_{j}(x,P^{k}_{i})
\right)_{j,i=1}^{n}$ of the corresponding systems are
nondegenerated and, consequently,
$a_{i(0,\dots,1,K-1)}=a_{i(0,\dots,K-1,1)}=0$, $i=1,\dots,n$.

Similarly, one can  show that all the leading coef\/f\/icients of the
operators $d_{k}$ are  zero. Thus, we have come to a
contradiction with the fact that $K$ is the leading order of
operators $d_{1},\dots,d_{n}$.
\end{proof}

To complete the proof of Theorem~\ref{T1} we need
\begin{lemma}\label{L3}
The modules $M_{\Gamma}$ and $N$ coincide.
\end{lemma}

\begin{proof} By $N(k)$ we denote  the following subset
\[
N(k)=\left\{\sum_{i=1}^{n}d_{i} \psi_{i}\:\vert\: d_{i}\in
\mathcal{D}_{n},\  {\rm ord}\,d_{i}\leq k-1\right\}.
\]
Since $\mathcal{D}_{n}$-module $N$ is free,
\[
 {\rm rank}_{\mathcal{O}}N(k)= n\; {\rm rank}_{\mathcal{O}}
 \left\{d \psi_{1}\, |\, m \; d\in \mathcal{D}_{n}, \ {\rm ord}\,d\leq k-1\right\}=
 n\, {\rm  rank}_{\mathcal{O}}\mathcal{D}_{n}=k  C_{k+n-1}^{n-1}.
\]
Consequently,
\[
 {\rm rank}_{\mathcal{O}}N(k)= {\rm rank}_{\mathcal{O}}M_{\Gamma}(k).
\]
Since there is  obvious inclusion $N(k) \subseteq
M_{\Gamma}(k)$, we obtain
\begin{gather*}
 M_{\Gamma}=N.\tag*{\qed}
\end{gather*}
\renewcommand{\qed}{}
\end{proof}
 Theorem~\ref{T1} is proven.

\subsection{Proof of subsidiary statements}\label{section3.3}

Choosing on ${\mathbb C}P^{1}\times{\mathbb C}P^{n-1}$ a suitable
coordinate system, we assume that $\Gamma$ is obtained by the
identif\/ication of hypersurfaces $(1:0)\times {\mathbb C}P^{n-1}$
and $(0:1)\times {\mathbb C}P^{n-1}$, i.e.~$a_2=b_1=0$. Also, we
can assume that $a_1=b_2=1$.
\begin{lemma}\label{L4}
Under condition \eqref{eq5} hypersurfaces $f=0$, $f_{1}=0$, $\dots$, $f_{n}=0$ do not have common points in ${\mathbb
C}P^{1}\times{\mathbb C}P^{n-1}$.
\end{lemma}

\begin{proof} Let us recall that $f$ and $f_1,\dots,f_n$ are the basis of
the form~(\ref{eq2}), satisfying the identity~(\ref{eq6}). We f\/ind such number
$k$ that $f_k$ satisf\/ies~(\ref{eq6}) for $c_k \neq 0$. Without the loss of
generality, we assume $c_{n} \neq 0$. Then the systems
\[
\left\{%
\begin{array}{l}
    f=0,\\
    f_1=0, \\
    \cdots\cdots\cdots\\
    f_{n-1}=0, \\
    f_n=0,
\end{array}%
\right. \qquad \mbox{and}
\qquad \left\{%
\begin{array}{l}
    f=0,\vspace{1mm}\\
   \displaystyle f_1-\frac{c_1}{c_n}f_n=0 ,\\
    \cdots\cdots\cdots\cdots\cdots\\
 \displaystyle   f_{n-1}-\frac{c_{n-1}}{c_n}f_n=0, \vspace{1mm}\\
 \displaystyle   \frac{1}{c_n}f_n=0,
\end{array}
\right.
\]
are equivalent. It is easy to verify that $f_k-\frac{c_k}{c_n}f_n$
satisfy (\ref{eq3}). Let $\tilde{f}_k=f_k-\frac{c_k}{c_n}f_n$ for $1
\leq k \leq n-1,$ and $\tilde{f}_n=\frac{1}{c_n}f_n$.

By the def\/inition of  $f$ and condition (\ref{eq3}), we have
\begin{gather*}
f(z_{1},z_{2},t)= f(z_{1},0,t)+f(0,z_{2},t)= z_{1} f(1,0,t)+z_{2}
f(0,1,t)\\
\phantom{f(z_{1},z_{2},t)}{}=
z_{1}  A  f(0,1,\mathcal{P}(t))+z_{2}f(0,1,t).
\end{gather*}
Since $f(z,t)$ is linear in last $n$ arguments, then the equality
$f=0$ can be rewritten in the following way:
\[
f(0,1,z_{1} A\mathcal{P}(t)+z_{2}t)=0.
\]
Likewise, we transform the equalities $\tilde{f}_{1}=0$, $\dots$,
$\tilde{f}_{n-1}=0$. We obtain the following system of linear
equations
\[
    f(0,1,v)=0,     \qquad
    \tilde{f}_{j}(0,1,v)=0, \qquad j=1,\dots,n-1,
\]
 $v=(z_{1}  A \mathcal{P}(t)+z_{2} t) \in {\mathbb C}^{n}$. By
linear independence (over ${\mathbb C}$) of the functions $f,\widetilde{f}_{1},\dots,\widetilde{f}_{n-1}$ this system has the
unique solution $v=0$. Thus, the equations $f=0$,
$f_{1}=0$, $\dots$, $f_{n-1}=0$ are equivalent to
\[
z_{1} A \mathcal{P}(t)+z_{2} t =0.
\]
The solutions of this system have the form:
\begin{gather*}
1)  \  z_1=z_2=0, \\
2)  \ t=0,\\
3) \ z_{2} /(A z_{1})=-\lambda_{j}, \qquad t=\mathbf{w}_{j}.
\end{gather*}
The solutions of the form 1) and 2) do not specify any point in
${\mathbb C}P^{1}\times{\mathbb C}P^{n-1}$. Taking into account
(\ref{eq4}),  case 3) gives $n$ dif\/ferent solutions. Now we add to this
system the last equality $\tilde{f}_{n}=0$. According to~(\ref{eq7}),
\begin{gather*}
\tilde{f}_{n}(z,t)= z_{1} A \tilde{f}_{n}(0,1,\mathcal{P}(t))+
z_{2}\tilde{f}_{n}(0,1,t)+f(1,0,t)=
\tilde{f}_{n}(0,1,z_{1} A \mathcal{P}(t)+z_{2}\,t)+f(1,0,t).
\end{gather*}
Taking into account (\ref{eq5}) and 3), we have
\[
\tilde{f}_{n}(z,t)=\tilde{f}_{n}(0,1,0,\dots,0)+f(1,0,\mathbf{w}_{j})=f(1,0,\mathbf{w}_{j})\neq
0.
\]
This reasoning ends the proving of Lemma~\ref{L4}.
\end{proof}

\begin{proof}[Proof of Lemma~\ref{L2}] We note that since  any function of the form (\ref{eq2}) is
a linear combination of the functions $f$ and $f_{1},\dots,f_{n}$,
then from the def\/inition of the module $M_{\Gamma}$ it follows
that a dif\/ferent choice of the forms $f_j$ corresponds to a
nondegenerated linear change of variables $x_1,\dots,x_n$. This
implies that if the statement of Theorem 1 is true for any f\/ixed
set $f_{1},\dots,f_{n}$ then it is true for any other one.
Therefore, we work with easy-to-use basis.

Let us choose the basis such that the points of the intersection
$P_{i}^{k}=(z_{i 1}^{k},z_{i 2}^{k},t_{i}^{k})$ of the
hypersurfaces $f(P)=0$ and $f_{j}(P)=0$ ($j=1,\dots,n$; $j\neq k$)
satisfy the following two conditions:
\begin{enumerate}\itemsep=0pt
\item[A)] f\/irst two coordinates $z_{i 1}^{k}$ and $z_{i 2}^{k}$ of the
points $P_i^k$ do not vanish;

\item[B)] the set of the vectors $(t_1^k,\dots,t_n^k) \in {\mathbb C}^n$ is
 linearly independent.
 \end{enumerate}

It always can be achieved. Indeed, we choose $n$ functions
$\tilde{f}_{j}$ of the form (\ref{eq2}), satisfying (\ref{eq3}), and such that
any $n-1$ from them and the function $f$ are linear independent. In
the course of proving of Lemma 4 it was shown that the points
of intersection of the hypersurfaces $f=0$, $\tilde{f}_{j}=0$,
$j=1,\dots,n$; $j\neq k$, satisfy the conditions~A) and~B). We take
one more function $\tilde{f}_{n+1}$ of the same form, satisfying~(\ref{eq6}) at $\tilde{c}_{n+1}=1$. Since the coordinates of the points of
intersections continuously depend on the coef\/f\/icients of the
functions, then by $f_{j}$ we take $f_{j}=\tilde{f}_{j}+c_{j}
\tilde{f}_{n+1}$, where $c_{j}$ are suf\/f\/iciently small. From the
proof of Lemma~\ref{L4} we obtain that all points $P^{k}_{j}$ are
dif\/ferent.

Let us write  matrix $(h_j(x,P_i^k))$ in a more convenient
form. For this, we note that the condi\-tion~(\ref{eq8}) can be written as
the condition on $h_{k}$
\[
h_{k}(x,a_{1},b_{1},t)-h_{k}(x,a_{2},b_{2},\mathcal{P}(t)) \Lambda
A e^{-cx}=0.
\]
Then the following equalities are true:
\begin{gather}
h_{j}(x,z_{1},z_{2},t)= h_{j}(x,z_{1},0,t)+h_{j}(x,0,z_{2},t)
=
z_{1} \Lambda A e^{-cx}  h_{j}(x,0,1,\mathcal{P}(t))
+z_{2} h_{j}(x,0,1,t)
\nonumber\\
\phantom{h_{j}(x,z_{1},z_{2},t)}{}  = h_{j}(x,0,1,z_{2} t+z_{1} \Lambda A
e^{-cx} \mathcal{P}(t)).\label{eq18}
\end{gather}
Let $\mathbf{v}^{k}_{i}=\left(z^{k}_{i 2} t^{k}_{i}+
z_{i 1}^{k}  e^{-cx} \mathcal{P}(t^{k}_{i})\right)$. Thus, the
nondegeneracy condition of the matrix $(h_{j}(x,P^k_{i}))$ is
written in the following way
\begin{equation}\label{eq19}
\det\left(%
\begin{array}{cccc}
  h_{1}(x,0,1,\mathbf{v}^k_{1}) & h_{2}(x,0,1,\mathbf{v}^k_{1}) & \dots & h_{n}(x,0,1,\mathbf{v}^k_{1}) \\
  \vdots                          & \vdots                          & \ddots & \vdots \\
  h_{1}(x,0,1,\mathbf{v}^k_{n}) & h_{2}(x,0,1,\mathbf{v}^k_{n}) & \dots & h_{n}(x,0,1,\mathbf{v}^k_{n}) \\
\end{array}%
\right) \neq 0.
\end{equation}
Since $h_{j}(x,P)$ are independent over $\mathcal{O}$, from (\ref{eq18}) we obtain that the functions $h_j(x,0,1,\cdot)$, as the functions
of the last $n$ arguments, are also independent over
$\mathcal{O}$. Then inequality (\ref{eq19}) is equivalent to the linear
independence of the vectors $\mathbf{v}^k_{i}$, since $h_i$ are the
linear forms of $\mathbf{v}^k_{i}$.

Let us show that the vectors $\mathbf{v}^k_{i}$ are linear independent.
Suppose that it is wrong, i.e.\ there are coef\/f\/icients
$\gamma_{i}$ (generally speaking dependent on $x$) such, that
$\sum \gamma_{i} \mathbf{v}^k_{i}=0$ or in more detail:
\[
\sum_{i=1}^{n}\left( \gamma_{i} z_{i 2} t_{i}+ \gamma_{i}
z_{i 1} \Lambda A e^{-cx} s_{i}\right)=0,
\]
where $s^k_{i}=\mathcal{P}(t^k_{i})$. The last equality can be
written in the matrix form
\begin{equation}\label{eq20}
\left(T+ \Lambda A  e^{-cx} S\right) \gamma=0,
\end{equation}
where $T$ and $S$ are matrices, composed from vectors
$z_{i 2} t^k_{i}$ and $z_{i 1} s^k_{i}$ respectively and
independent from $x$, $\gamma=\left(\gamma_{1},\dots, \gamma_{n}
\right)^{\top}$. Since the matrices $T$ and $S$ are nondegenerated,
then there are not more then $n$ values $\mu_{j}$ such, that
$\det \left(T+\mu_{j} S\right) = 0$ and $\mu_{j}$ also do not
depend on $x$. Consequently, for $\Lambda A   e^{-cx} \neq
\mu_{j}$ system (\ref{eq20}) has only one solution $\gamma=0$, i.e.,
$\mathbf{v}^k_{i}$ are linear independent for almost every $x$
and, consequently, the determinant $|h_{j}(x,P^k_{i})|$ does not
vanish identically in~$x$. Lemma~\ref{L2} is proven.
\end{proof}

\section[Proof of Theorem 2]{Proof of Theorem~\ref{T2}}\label{section4}

Let us choose in $M_{\Omega}(1)$ two independent over $\mathcal{O}$
functions $\varphi_1$ and $\varphi_2$
\[
 \varphi_i=\frac{\tilde{h}_i(x,y,P)}{g(P)}
 \exp\left(xG_1(P)+yG_2(P)\right),
\]
where
\[
 \tilde{h}_i(x,y,P)=k_i(x,y)z_1z_2+l_i(x,y)z_1w_2+m_i(x,y)w_1z_2+n_i(x,y)w_1w_2,\qquad
 i=1,2.
\]
The functions $\tilde{h}_i$ satisfy the identity
\[
 \tilde{h}_i(x,y,1:0,t_1:t_2)-\tilde{h}_i(x,y,t_1:t_2,0:1)\Lambda B
 \exp(-xc_1-yc_2)=0.
\]

By $P_i$ and $Q_i$ we denote  the points of intersection of the
curves, def\/ined by the equations $g_i(P)=0$ and $g(P)=0$, $i=1,2$.
By the inf\/initesimal variations $c_1$ and $c_2$, one can obtain
that these points are pairwise dif\/ferent.

By means of direct check, we can ascertain that the determinant
\[
\det \left(%
\begin{array}{cc}
  \tilde{h}_{1}(x,y,P_i) & \tilde{h}_{2}(x,y,P_i) \\
  \tilde{h}_{1}(x,y,Q_i) & \tilde{h}_{2}(x,y,Q_i) \\
\end{array}%
\right)
\]
up to multiplication by nonvanishing function from $\mathcal{O}$,
is equal to
\begin{equation}\label{eq21}
\frac{\sqrt{(\gamma\delta_i-\gamma_i\delta)^2+2 \gamma \delta c_i
(\gamma\delta_i-\gamma_i\delta)+ \delta^2(\gamma-2 B \delta)^2
c_i^2}(\Lambda e^{-c_1 x-c_2
y}-1)^2}{\gamma\delta_i-\gamma_i\delta}.
\end{equation}
It is obvious that the inequality
\begin{equation}\label{eq22}
\det \left(%
\begin{array}{cc}
  \tilde{h}_{1}(x,y,P_i) & \tilde{h}_{2}(x,y,P_i) \\
  \tilde{h}_{1}(x,y,Q_i) & \tilde{h}_{2}(x,y,Q_i)
\end{array}%
\right) \neq 0
\end{equation}
is fulf\/illed for almost all $x$ and $y$.

Further the proof of Theorem~\ref{T2} verbatim repeats the proof of
Theorem~\ref{T1}.

\section{Examples}\label{section5}

In this section we demonstrate the examples of the commuting
dif\/ferential operators and their common eigenvector-functions
def\/ining the basis in the free BA-modules.

\subsection[Commuting operators corresponding to the variety $\Gamma$]{Commuting operators corresponding to the variety $\boldsymbol{\Gamma}$}

Let us consider the case $n=2$. As a spectral variety we take
\[
\Gamma={\mathbb C}P^{1}\times{\mathbb C}P^{1} /
\{(1:0,t_{1}:t_{2})\sim(0:1,t_{2}:t_{1})\},
\]
i.e.\ in terms of Section~\ref{section2}  $ p_{1}=(1:0)$, $p_{2}=(0:1)$,
$\mathcal{P}(t_{1},t_{2})=(t_{2},t_{1})$.

We introduce three functions
\begin{gather*}
f(P)=-z_{1} t_{1}-z_{2} t_{2},
\qquad
f_{1}(P)=z_{1} t_{2}+z_{2} ( t_{1}-i t_{2}),
\qquad
f_{2}(P)=-z_{1} t_{2}-z_{2} ( t_{1}+i t_{2}).
\end{gather*}
Via the direct check, we can ascertain that  $f(P)$
satisf\/ies the condition (\ref{eq3}) for $A=1$, and~$f_{1}$,~$f_{2}$ satisfy the condition (\ref{eq6}) for $c_{1}=c_{2}=-i$.

Let us choose in $\mathcal{D}$-module $M_{\Gamma}$ the basis
\begin{gather*}
\psi_{1}(x,y,P) = \frac{z_{1}t_{1}+e^{-i(x+y)}z_{2}t_{2}}{f(P)}
\exp\left(\frac{f_{1}}{f}x+\frac{f_{2}}{f}y\right),
\\
\psi_{2}(x,y,P) = \frac{z_{1}t_{2}+e^{-i(x+y)}z_{2}t_{1}}{f(P)}
\exp\left(\frac{f_{1}}{f}x+\frac{f_{2}}{f}y\right).
\end{gather*}
We consider the following meromorphic functions on $\Gamma$ with the
poles on the curve $f(P)=0$
\begin{gather*}
\lambda_{1}=\frac{2(z_{1}t_2+z_{2}t_{1})}{f(P)},\qquad
\lambda_{2}=\frac{i z_{1}z_{2}(-t_{1}^{2}+t^{2}_{2})}{f(P)^{2}},\qquad
\lambda_{3}=\frac{z_{1}^{2}t_{2}^{2}+3
z_{1}z_{2}t_{1}t_{2}+z_{2}^{2}t_{1}^{2}}{f(P)^{2}}.
\end{gather*}
Pairwise commuting operators, corresponding to these functions
have the forms
\begin{gather*}
D(\lambda_{1})=\left(%
\begin{array}{cc}
  \partial_{x}-\partial_y & 0              \\
  0              & \partial_{x}-\partial_y
\end{array}
\right),
\\
D(\lambda_{2})=\left(%
\begin{array}{cc}
\displaystyle   \frac{1}{4}\big(\partial_{y}^{2}-\partial_{x}^{2}\big)+\frac{1}{2}\cot\left(\frac{x+y}{2}\right)(\partial_{x}-\partial_{y}) &
 \displaystyle \cot\left(\frac{x+y}{2}\right)-\frac{1}{2}(\partial_{x}+\partial_{y})\vspace{2mm}\\
\displaystyle   -\frac{1}{2}(\partial_{x}+\partial_{y}) &
\displaystyle  \frac{1}{4}\big(\partial_y^2-\partial_{x}^2\big)
\end{array}
\right).
\end{gather*}
Operator, corresponding to the function $\lambda_3$ has the form
\begin{gather*}
[D(\lambda_3)]_{11} =
\frac{1}{2}\partial_{x}^{2}+\frac{1}{2}\partial_{y}^{2}-
\frac{1}{2}\cot\left(\frac{x+y}{2}\right)(\partial_{x}+\partial_{y}),
\\
 [D(\lambda_3)]_{12} =0,\qquad  [D(\lambda_3)]_{21}=
 \frac{1}{4 \sin^{2}\big(\frac{x+y}{2}\big)}(\partial_{x}-\partial_{y}),
\\
 [D(\lambda_3)]_{22}=
 \frac{1}{2 \sin^{2}\big(\frac{x+y}{2}\big)}-\frac{1}{2}\cot\left(\frac{x+y}{2}\right)(\partial_{x}+\partial_{y})+
 \frac{1}{2}\partial_{x}^{2}+\frac{1}{2}\partial_{y}^{2}.
\end{gather*}

\subsection[Commuting operators corresponding to the variety $\Omega$]{Commuting operators
corresponding to the variety $\boldsymbol{\Omega}$}

Let us consider three functions
\begin{gather*}
 g(P)=z_1w_1+z_1w_2+z_2w_2,
\\
 g_1(P)=z_1w_1+2 z_2w_1-z_2w_2,
\\
 g_2(P)=-z_1w_1+2z_2w_1+z_2w_2.
\end{gather*}
By the direct check we can ascertain that  $g(P)$
satisf\/ies the identity (\ref{eq11}) for $B=1$,   $g_1(P)$
and $g_2(P)$ satisfy to the identity (\ref{eq12}) for $c_1=1$ and $c_2=-1$
respectively.

The curves $g(P)=0$ and $g_1(P)=0$ are intersected in the points
\begin{gather*}
 P_1=\left(-2-\sqrt{2}:1,-\frac{1}{\sqrt{2}}:1\right),\qquad
 Q_1=\left(-2+\sqrt{2}:1,\frac{1}{\sqrt{2}}:1\right),
\end{gather*}
and the curves  $g(P)=0$ and $g_2(P)=0$ are intersected in the
points
\[
 P_2=\left(-\sqrt{2}:1,-1+\frac{1}{\sqrt{2}}:1\right),\qquad
 Q_2=\left(\sqrt{2}:1,-1-\frac{1}{\sqrt{2}}:1\right).
\]
Chose the basis $\psi_1$, $\psi_2$ in $\mathcal{D}$-module
$M_{\Omega}$
\begin{gather*}
 \psi_1=\frac{z_2w_1}{g(P)}\exp\left(xG_1(P)+yG_2(P)\right),
\\
 \psi_2=\frac{z_1w_1e^{y-x}+z_1w_2+z_2w_2e^{x-y}}{g(P)}\exp\left(xG_1(P)+yG_2(P)\right).
\end{gather*}
Then
\begin{gather*}
 \frac{{\tilde{h}_1}(P_1,x,y)}{{\tilde{h}_2}(P_1,x,y)}=-\frac{e^{x+y}}{\sqrt{2}(e^y-e^x)(-e^x+(1+\sqrt{2})e^y)},
\\
 \frac{\tilde{h}_1(Q_1,x,y)}{\tilde{h}_2(Q_1,x,y)}=-\frac{e^{x+y}}{\sqrt{2}(e^y-e^x)(e^x+(-1+\sqrt{2})e^y)},
\end{gather*}
thus inequality (\ref{eq22}) is fulf\/illed.

Four  simplest meromorphic functions on $\Omega$ with the poles
on the curve $g(P)=0$ have the form
\[
  \lambda_1=\frac{z_2w_1}{g(P)},\qquad
  \lambda_2=\frac{z_1z_2w_1^2}{g(P)^2},  \qquad
  \lambda_3=\frac{z_1w_1z_2w_2}{g(P)^2},\qquad
  \lambda_4=\frac{z_1z_2w_2^2+z_1^2w_1w_2}{g(P)^2}.
\]
Pairwise commutating operators, corresponding to these functions
have the form
\begin{gather*}
 D(\lambda_1)=
  \left(
  \begin{array}{cc}
    \frac{1}{4}(\partial_x+\partial_y) &  0\vspace{1mm}\\
   0   &  \frac{1}{4}(\partial_x+\partial_y)
  \end{array}\right),
\\[1mm]
 [D(\lambda_2)]_{11}=\frac{e^x}{8(e^x-e^y)}\big(\partial_x^2-\partial_y^2\big)-
 \frac{e^{x+y}}{4(e^x-e^y)^2}(\partial_x+\partial_y),
\\[1mm]
 [D(\lambda_2)]_{12}=\frac{e^{x+y}}{16(e^x-e^y)^2}(\partial_x+\partial_y)^2,
\\
 [D(\lambda_2)]_{21}=\frac{1}{8}(e^{y-x}-e^{x-y}-2)\partial_x^2+
 \frac{1}{8}(e^{x-y}-e^{y-x}-2)\partial_y^2+\frac{1}{2}\partial_x\partial_y
\\[1mm]
\phantom{[D(\lambda_2)]_{21}=}{}+
 \frac{e^x+e^{2x-y}+5e^y-e^{2y-x}}{4(e^x-e^y)}\partial_x+
 \frac{3e^x-e^{2x-y}+3e^y+e^{2y-x}}{4(e^y-e^x)}\partial_y-
 \frac{e^y(2e^x+e^y)}{(e^x-e^y)^2},
\\[1mm]
 [D(\lambda_2)]_{22}=\frac{e^x}{8(e^y-e^x)}\partial_x^2-\frac{1}{4}\partial_x\partial_y+
 \frac{e^x-2e^y}{8(e^y-e^x)}\partial_y^2+
 \frac{e^y(2e^x+e^y)}{8(e^x-e^y)^2}(\partial_x+\partial_y).
\end{gather*}
The operator, corresponding to the function $\lambda_3$ has the
form
\begin{gather*}
 [D(\lambda_3)]_{11}=\frac{(e^x+e^y)}{8(e^y-e^x)}\big(\partial_x^2-\partial_y^2\big)+
 \frac{(e^{2x}+e^{2y})}{4(e^y-e^x)^2}(\partial_x+\partial_y),
\\[1mm]
 [D(\lambda_3)]_{12}=\frac{e^{x+y}}{8(e^y-e^x)^2}(\partial_x-\partial_y)^2,
\\[1mm]
 [D(\lambda_3)]_{21}=\frac{1}{4}(2+e^{x-y}-e^{y-x})\partial_x^2-
 \partial_x\partial_y+\frac{1}{4}(2-e^{x-y}+e^{y-x})\partial_y^2
\\[1mm]
 \phantom{[D(\lambda_3)]_{21}=}{}+
 \frac{2e^x+e^{2x-y}+4e^y-e^{2y-x}}{2(e^y-e^x)}\partial_x+
 \frac{4e^x-e^{2x-y}+2e^y+e^{2y-x}}{2(e^x-e^y)}\partial_y\\
 \phantom{[D(\lambda_3)]_{21}=}{}
 +
 \frac{e^{2x}+e^{2y}+4e^{x+y}}{(e^x-e^y)^2},
\\[1mm]
 [D(\lambda_3)]_{22}=\frac{3e^x-e^y}{8(e^x-e^y)}\partial_x^2+
 \frac{1}{2}\partial_x\partial_y+\frac{e^x-3e^y}{8(e^x-e^y)}\partial_y^2-
 \frac{3e^{x+y}}{2(e^y-e^x)^2}(\partial_x+\partial_y).
\end{gather*}
Operator corresponding to the function $\lambda_4$ has the form
\begin{gather*}
 [D(\lambda_4)]_{11}=\frac{e^x+3e^y}{4(e^x-e^y)}\partial_x^2+
 \frac{1}{2}\partial_x\partial_y-\frac{3e^x+e^y}{4(e^x-e^y)}\partial_y^2
\\[1mm]
\phantom{[D(\lambda_4)]_{11}=}{}-
 \frac{e^{2x}+3e^{2y}}{2(e^x-e^y)^2}\partial_x-\frac{3e^{2x}+e^{2y}}{2(e^x-e^y)^2}\partial_y-
 \frac{2e^{x+y}}{(e^y-e^x)^2},
\\[1mm]
 [D(\lambda_4)]_{12}=\frac{e^{x+y}}{2(e^y-e^x)^2}((\partial_x+\partial_y)^2+\partial_x+\partial_y),
\\[1mm]
 [D(\lambda_4)]_{21}=(e^{y-x}-e^{x-y}-2)\partial_x^2+4\partial_x\partial_y+(e^{x-y}-e^{y-x}-2)\partial_y^2
\\[1mm]
\phantom{[D(\lambda_4)]_{21}=}{}+
 \frac{5e^x+e^{2x-y}+9e^y-3e^{2y-x}}{e^x-e^y}\partial_x+
 \frac{9e^x-3e^{2x-y}+5e^y+e^{2y-x}}{e^x-e^y}\partial_y
\\[1mm]
\phantom{[D(\lambda_4)]_{21}=}{}+
 \frac{2e^{-x-y}(e^{4x}+e^{4y}-6e^{2(x+y)}-4e^{3x+y}-4e^{x+3y})}{(e^x-e^y)^2},
\\[1mm]
 [D(\lambda_4)]_{22}=\frac{7e^x-3e^y}{4(e^y-e^x)}\partial_x^2-
 \frac{3}{2}\partial_x\partial_y+\frac{3e^x-7e^y}{4(e^y-e^x)}\partial_y^2
\\[1mm]
\phantom{[D(\lambda_4)]_{22}=}{}-
\frac{e^{2x}+3e^{2y}-16e^{x+y}}{2(e^x-e^y)^2}\partial_x-
 \frac{3e^{2x}+e^{2y}-16e^{x+y}}{2(e^x-e^y)^2}\partial_y.
\end{gather*}

\appendix

\section[Structures of algebraic varieties on $\Gamma$ and $\Omega$]{Structures of algebraic varieties on $\boldsymbol{\Gamma}$ and $\boldsymbol{\Omega}$}\label{appendixA}

We show that on $\Gamma$ and $\Omega$ structures of algebraic varieties can be
introduced. For this in the f\/irst case we
construct a smooth morphism from ${\mathbb
C}P^{1}\times{\mathbb C}P^{n-1}$ to ${\mathbb C}P^{2}\times{\mathbb
C}P^{2n-1}$, and in the second case from ${\mathbb
C}P^1\times{\mathbb C}P^1$ to ${\mathbb C}P^{11}$. The
 morphisms are injective everywhere except
gluing hypersurfaces. The images of the morphisms are algebraic
varieties which def\/ine required structures on $\Gamma$ and $\Omega$.

\subsection[Variety $\Gamma$]{Variety $\boldsymbol{\Gamma}$}

By choosing the convenient coordinate system on ${\mathbb C}P^1$ we
can assume that
\[
\Gamma= {\mathbb C}P^{1}\times{\mathbb
C}P^{n-1}/\{(1:0,t)\sim(0:1,\mathcal{P}(t))\}.
\]
We consider the mapping
\[
\varphi_1: \ {\mathbb C}P^{1}\times{\mathbb C}P^{n-1} \rightarrow
{\mathbb C}P^{2}\times{\mathbb C}P^{2n-1},
\]
def\/ined by the formula:
\[
\varphi_1(z,t)=(u,v),
\]
where $u=(u_1:u_2:u_3)$,
\begin{gather*}
u_{1}=z_{1}^{2} z_{2},\qquad u_{2}=z_{1} z_{2}^{2}, \qquad u_{3}=z_{1}^{3} +
z_{2}^{3},
\qquad
v=\big(z_{1}^{2} t + z_{2}^{2} \mathcal{P}^{-1}(t):z_{1}z_{2} t\big).
\end{gather*}
Here $v = (\xi_1:\dots:\xi_n:\eta_1:\dots:\eta_n)$, $\xi_j=z_12
t_j+z_22 r_j$ and $\eta_j=z_1 z_2 t_j$, where $r_j$ is $j$-th
coordinate $\mathcal{P}^{-1}(t)$, $j=1,\dots,n$.

\begin{lemma} \label{L5}
The mapping $\varphi_1$ is correctly defined on $\Gamma$ and
is the embedding of $\Gamma$. The image of~$\varphi_1$ is defined
by the equations:
\begin{gather}\label{eq23}
u_{1}^{3}+u_{2}^{3} = u_{1}u_{2}u_{3},
\\
\label{eq24}
u_{1}^{2} \eta + u_{2}^{2} \mathcal{P}^{-1}(\eta)=u_{1}u_{2} \xi,
\\
\label{eq25}
u_{3}\eta+u_{1} \mathcal{P}(\eta) +
u_{2}\mathcal{P}^{-1}(\eta)=u_{1} \xi +u_{2} \mathcal{P}(\xi).
\end{gather}
\end{lemma}
\begin{proof} Let us show that the image of the point $(z,t) \in {\mathbb
C}P^{1}\times{\mathbb C}P^{n-1}$ satisf\/ies  the equations
(\ref{eq23})--(\ref{eq25}). The equalities (\ref{eq23}) and (\ref{eq24}) obviously follow from the
def\/inition of $\varphi_1$. The equali\-ty~(\ref{eq25}) for $z_1 \neq 0$ and
$z_2 \neq 0$ is the corollary of (\ref{eq23}) and (\ref{eq24}). Indeed, from (\ref{eq24})
we obtain
\[
u_{1}^{2} \mathcal{P}(\eta) + u_{2}^{2} \eta=u_{1}u_{2}
\mathcal{P}(\xi).
\]
Let us multiply the obtained equality by $u_2$, equality (\ref{eq24}) by $u_1$
and take a sum. Dividing the result by $u_1u_2 \neq 0$, we obtain
(\ref{eq25}). If $z_1$ or $z_2$ are equal to zero, then $\eta=0$, and,
consequently, the left and right parts of (\ref{eq25}) vanish.

We show that for any point $(u,v) \in {\mathbb
C}P^{2}\times{\mathbb C}P^{2n-1}$, satisfying (\ref{eq23})--(\ref{eq25}), the inverse image can be
found. Note that from (\ref{eq23}) it follows that
$u_{1}$ and $u_{2}$ can vanish only simultaneously.

If $u_{1}=u_{2}=0$, then from (\ref{eq25}) we obtain $\eta=0$ and the
inverse image can be either  point $(1:0,\xi)$, or  point
$(0:1,\mathcal{P}(\xi))$, which are identif\/ied in $\Gamma$.

Now we consider the case $u_{1}\neq 0$ and $u_{2}\neq 0$. From
(\ref{eq24}) it follows that $\xi \neq 0$. In this case as the inverse
image of the point we can take  point $(u_{1}:u_{2},\eta)$.
 By easy calculations it can be checked
that $\varphi_1(u_{1}:u_{2},\eta)= (u_{1}:u_{2}:u_{3},\xi:\eta)$.
If any other point $B=(a:b,s_{1}:\dots:s_{n})$ is the inverse
image of $(u_{1}:u_{2}:u_{3},\xi:\eta)$, then $a:b=u_{1}:u_{2}$
and, consequently, $a\neq 0$, $b \neq 0$, $\eta_{j}=abs_{j}$,
i.e.\ $B = (u_{1}:u_{2},\eta)$, which shows injectivity of the
mapping $\varphi_1$ on $\Gamma$.

By direct calculations one can show that the dif\/ferential of
the mapping $\varphi_1$ is nondegene\-rated.
\end{proof}

\subsection[Variety $\Omega$]{Variety $\boldsymbol{\Omega}$}

Let us consider the mapping
\[
\varphi_2 : \ {\mathbb C}P^1 \times {\mathbb C}P^1 \rightarrow
{\mathbb C}P^{11},
\]
def\/ined by the formula:
\[
\varphi_2(z_1:z_2, w_1:w_2) = (u_1:\dots:u_{12}),
\]
where
%\begin{alignat*}{4}
%& u_1 = z_1^3 \big(w_1^3 + w_2^3\big) + (z_2 w_2)^3, \qquad && u_2 = z_1^3 w_1^2 w_2 + z_1^2 z_2 w_2^3,\qquad &&
% u_3 = z_1^3 w_1 w_2^2 + z_1 z_2^2 w_2^3, & \\
%& u_4 = z_1^2 z_2 w_1^3,  && u_5 = z_1^2 z_2 w_1^2 w_2, && u_6 = z_1^2 z_2 w_1 w_2^2, & \\
%& u_7 = z_1 z_2^2 w_1^3,  && u_8 = z_1 z_2^2 w_1^2 w_2, && u_9 = z_1 z_2^2 w_1 w_2^2, & \\
%& u_{10} = z_2^3 w_1^3, && u_{11} = z_2^3 w_1^2 w_2, && u_{12} = z_2^3 w_1 w_2^2.
%\end{alignat*}
\begin{gather*}
  u_1 = z_1^3 \big(w_1^3 + w_2^3\big) + (z_2 w_2)^3, \qquad   u_2 = z_1^3 w_1^2 w_2 + z_1^2 z_2 w_2^3,\qquad
 u_3 = z_1^3 w_1 w_2^2 + z_1 z_2^2 w_2^3,   \\
  u_4 = z_1^2 z_2 w_1^3,\!\!  \qquad u_5 = z_1^2 z_2 w_1^2 w_2,\!\! \qquad u_6 = z_1^2 z_2 w_1 w_2^2, \!\!\qquad
  u_7 = z_1 z_2^2 w_1^3, \!\! \qquad u_8 = z_1 z_2^2 w_1^2 w_2, \\
   u_9 = z_1 z_2^2 w_1 w_2^2, \qquad
  u_{10} = z_2^3 w_1^3, \qquad u_{11} = z_2^3 w_1^2 w_2, \qquad u_{12} = z_2^3 w_1 w_2^2.
\end{gather*}
We  can easily ascertain that $u_j$ do not vanish simultaneously
(for example, $u_1$, $u_2$, $u_3$ and $u_{10}$ equal zero if and
only if $z_1=z_2=0$ or $w_1=w_2=0$).

\begin{lemma}\label{L6}
The mapping $\varphi_2$ is correctly defined on $\Omega$ and is the
embedding of $\Omega$.
\end{lemma}

\begin{proof} We show that the mapping $\varphi_2$ identif\/ies only
points $(1:0,t_1:t_2)$ and $(t_1:t_2,0:1)$ on ${\mathbb C}P^1
\times {\mathbb C}P^1$.

If $u_{10} \neq 0$, then from the def\/inition of $\varphi_2$ it
follows that the inverse image has the form $(z_1:z_2,w_1:w_2) =
(u_7:u_{10},u_{10}:u_{11})$.

If $u_{10} = 0$, then $u_k = 0$ for $4\le k \le 12$, $u_2$ and
$u_3$ vanish simultaneously. Two cases are possible:
%\begin{enumerate}\itemsep=0pt

\textbf{a)} $u_2 \ne 0$ and $u_3 \neq 0$, then inverse image is one of the two
       points $(1:0,u_2:u_3)$ or $(u_2:u_3,0:1)$, which are identif\/ied in $\Omega$;

\textbf{b)} $u_2 = u_3 = 0$, then $u_1 \neq 0$ and inverse image is one of the three points
       $(1:0,1:0)$, $(1:0,0:1)$ or $(0:1,0:1)$, which are also identif\/ied in $\Omega$.

By  direct calculations one can show that the dif\/ferential of
the mapping $\varphi_2$ is nondegene\-rated.
\end{proof}

\subsection*{Acknowledgements}

The work is supported by the Program of
Russian Academy of Sciences ``Fundamental Problems of Nonlinear
Dynamics''.
 The second author (A.E.M) is thankful to
Atsushi Nakayashiki for the invitations in Kyushu University,
useful discussions of our results. The second author is also
grateful to Grant-in-Aid for Scientif\/ic Research (B) 17340048 for
f\/inancial support of the visits in Kyushu University.

\pdfbookmark[1]{References}{ref}
\LastPageEnding

\end{document}